\documentclass[10pt,twoside]{article}
\usepackage{amssymb}
\usepackage{psfig}

\newtheorem{lemma}{Lemma}

\newtheorem{definition}[lemma]{Definition}

\newtheorem{conjecture}[lemma]{Conjecture}

\begin{document}

\renewcommand{\thefootnote}{\fnsymbol{footnote}}
\setcounter{page}{1}

\bigskip
\bigskip

\large
\centerline{\bf KNOT POLYNOMIALS: MYTHS AND REALITY}
\normalsize

\bigskip
\bigskip
\bigskip

\bigskip
\centerline{\bf Slavik Jablan, Ljiljana Radovi\' c$^*$}

\bigskip

\centerline{\footnotesize\it The Mathematical Institute, Knez
Mihailova 36,}\centerline{\footnotesize\it P.O.Box 367, 11001
Belgrade, Serbia} \centerline{\footnotesize\it sjablan@gmail.com}

\bigskip

\centerline{\footnotesize\it Faculty of Mechanical
Engineering$^{*}$, A.~Medvedeva 14,}
\centerline{\footnotesize\it 18 000 Ni\v s, Serbia} \centerline{\footnotesize\it ljradovic@gmail.com}

\normalsize


\renewcommand{\thefootnote}{\roman{footnote}}

\begin{abstract}
This article provides an overview of relative strengths of
polynomial invariants of knots and links, such as the Alexander,
Jones, Homflypt, and Kaufman two-variable polynomial, Khovanov
homology, factorizability of the polynomials, and knot primeness detection.
\end{abstract}

\section{Introduction}

In some sources the end of the 19$^{th}$ century is called the "dark
age of the knot theory", because knots and links ($KL$s) are
recognized "by hand" or some other "non-exact methods". However,
first knot tables were created at this time by P.G. Tait, T.P.
Kirkman and C.N. Little, after more than five years of a hard work.
In knot tabulation, almost nothing important happened almost a
century, until the computer derivation of knot and link tables by M.
Thistlethwaite and his collaborators [1], and now computations have
reached the limit even with the use of supercomputers. Let us give
the overview of the polynomial invariants we have at hand.

The first knot polynomial introduced by J.W. Alexander was used by
K. Reidemeister in his book {\it Knotentheorie} in 1932 to
distinguish knots up to $n=9$ crossings. A new series of invariants,
beginning with the Jones polynomial, is recently extended by using
categorifications to the more powerful invariants.

Appearance of every knot invariant is usually connected with the
progress in different fields of mathematics (e.g., the Alexander
polynomial and Fox calculus, Khovanov homology and categorifications
in different fields of algebra) and its connections with other
sciences, in particular with physics (e.g., the Jones polynomial and
its relation with the Potts model). In this paper we will not
discuss the impact of knot polynomials to the development of
mathematics or other fields of science, but only their ability to
distinguish different $KL$s.

One of the first things we learn in knot theory is the computation
of polynomial knot invariants, mostly those that can be computed by
using skein relations. After learning that the Alexander polynomial
is not able to distinguish a left trefoil from the right, that it
cannot recognize unknot, that the Jones polynomial can distinguish
left and right trefoil and (maybe) recognizes unknots, we believe
that we have in our hand a very powerful tool for knot recognition,
despite of the fact (usually illustrated by a few standard examples)
that for every polynomial invariant exist $KL$s (not only mutants)
that it cannot distinguish.

For all computations we used the program {\it LinKnot} [2], combined
with the programs [1,3,4].

\section{Distinction of knots and links by polynomial \\ invariants}

In order to compare different polynomial invariants and their
ability to distinguish different $KL$s we computed different $KL$
polynomials for all $KL$s up to $n=12$ crossings and the number of
$KL$s sharing the same polynomial with some other $KL$. Because
there are 4684 alternating $KL$s with $n\le 12$ crossings,
consisting of 1851 knots and 2833 links, and 3993 non-alternating
$KL$s with $n\le 12$ crossings consisting of 1126 knots and 2867
links, i.e., 8677 $KL$s in total, we believe that this is a large
enough sample from which we can make some conclusions.

In the following tables is given the name of the corresponding
polynomial, number of knots sharing the same polynomial with some
other knot, their percent among all knots, the same results for
links, and the total number and percent of $KL$s that cannot be
distinguished by the corresponding polynomial. The Table 1 contains
the data about alternating, Table 2 about non-alternating $KL$s, the
Table 3 is the sum of Table 1 and Table 2, and Table 4 shows the
results of computations for all $KL$s, where alternating $KL$s are
not separated from non-alternating ones.

\bigskip
\medskip

{\bf Table 1}

\medskip

\begin{center}
\begin{tabular}{|c|c|c|c|c|c|c|} \hline
 \bf{Alternating} &Knots & &Links & &Total & \\ \hline
 Alexander & 846 & 46\%& 1732& 61\%& 2578& 55\% \\ \hline
 Jones  & 601 & 32\%& 672 & 24\%& 1273& 27\%\\ \hline
 Khovanov & 599 & 32\%& 406 & 14\%& 1005 & 21\%\\ \hline
 Homflypt & 274& 15\%&  285& 10\%& 559 & 12\%\\ \hline
 Kauffman&  93 & 5\%& 243& 9\%& 336& 7\%\\ \hline
 \end{tabular}
 \end{center}
\medskip
\medskip

In our computation are not included some very powerful $KL$ invariants: colored
Jones polynomials and Links-Gould invariant, which cannot be
computed for so large amount of $KL$s in a reasonable
time. In the recognition of $KL$s, odd Khovanov homology gives the same results as the Khovanov homology\footnote{The authors are thankful to Krzystof Putyra for noticing the errors in the computations of Khovanov and odd Khovanov homology, that appeared in the first version of this paper.}.

\bigskip

{\bf Table 2}

\medskip

\begin{center}
\begin{tabular}{|c|c|c|c|c|c|c|} \hline
\bf{Non-alternating} &Knots & &Links & &Total & \\ \hline
 Alexander & 697& 62\%& 2123& 74\%& 2820& 71\%\\ \hline
 Jones  & 459&  41\%& 797& 28\%& 1256& 31\%\\ \hline
 Khovanov  & 398& 35\%& 459& 16\%& 857 & 21\%\\ \hline
 Homflypt  & 254& 23\%& 400& 14\%& 654& 16\%\\ \hline
 Kauffman  & 146& 13\%& 327& 11\%& 473& 12\%\\ \hline
\end{tabular}
\end{center}
\medskip

{\bf Table 3}

\medskip

\begin{center}
\begin{tabular}{|c|c|c|c|c|c|c|} \hline
\bf{Sum} &Knots & &Links & &Total & \\ \hline
 Alexander & 1543& 52\%& 3855& 68\%& 5398& 62\%\\ \hline
 Jones  & 1060&  36\%& 1469& 26\%& 2529& 29\%\\ \hline
 Khovanov  & 997 & 33\%& 865 & 15\%& 1862& 21\%\\ \hline
 Homflypt  & 528& 18\%& 685& 12\%& 1213& 14\%\\ \hline
 Kauffman  & 239& 8\%& 570& 10\%& 809& 9\%\\ \hline
\end{tabular}
\end{center}
\medskip

{\bf Table 4}

\medskip

\begin{center}
\begin{tabular}{|c|c|c|c|c|c|c|} \hline
\bf{All} &Knots & &Links & &Total & \\ \hline
 Alexander & 1832& 62\%& 4169& 73\%& 6001& 69\%\\ \hline
 Jones  & 1213&  41\%& 1565& 27\%& 2778& 32\%\\ \hline
 Khovanov  & 1117 & 38\%& 921 & 16\%& 2038 & 23\%\\ \hline
 Homflypt  & 600& 20\%& 707& 12\%& 1307& 15\%\\ \hline
 Kauffman  & 239& 8\%& 570& 10\%& 809& 9\%\\ \hline
\end{tabular}
\end{center}

\medskip

\begin{figure}
\centerline{\psfig{file=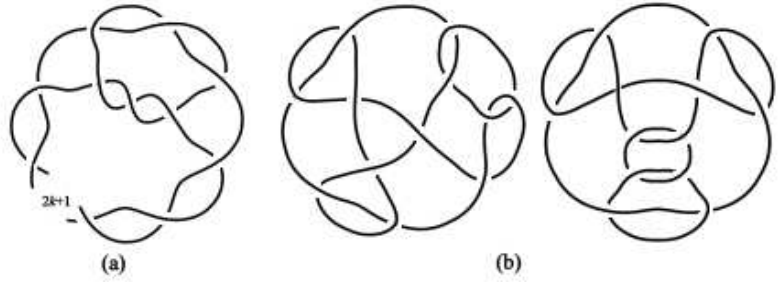,width=4.20in}} \vspace*{8pt}
\caption{(a) Knot family $(2k+1),3,-3$; (b) 2-component link
$(2\,1,2\,1)\,1\,(2,2+)$; (c) 4-component link
$(2,2,2)\,(2\,1,2\,1)$.} \label{f1}
\end{figure}

\begin{definition}
For a link $L$ given in an unreduced \footnote{The Conway notation
is called unreduced if in symbols of polyhedral links elementary
tangles 1 in single vertices are not omitted.} Conway notation
$C(L)$, let $S$ denote a set of numbers in the Conway symbol,
excluding numbers denoting basic polyhedra and zeros (marking the
position of tangles in the vertices of polyhedra), and $S_f$ the set
obtained by substituting every positive number from $S$ different
from $1$ by $2$, and every negative number from $S$ different from
$-1$ by $-2$. For $C(L)$ and an arbitrary (non-empty) subset $\tilde
S$ of $S$ the family $F_{\tilde S}(L)$ of knots or links derived
from $L$ is constructed by substituting each $a \in S_f$, $a\neq 1$,
by $sgn(a) (|a|+k_a)$ for $k_a \in N$.
\end{definition}

If $k_a$ is an even number ($k_a\in N$), the number of components is
preserved inside a family, i.e., we obtain families of knots or
links with the same number of components.

For the Alexander polynomial, there are even families of knots that
can not be distinguished one from another. For example, for all
knots of the family of nonalternating pretzel knots $(2k+1),3,-3$
(Fig. 1a), the Alexander polynomial is $2-5x+2x^2$.

All polynomials distinguish knots from links, but Alexander
polynomial cannot distinguish links according to the number of
components, and all the other polynomials distinguish them. For
example, 2-component link $(2\,1,2\,1)\,1\,(2,2+)$ with $n=12$
crossings and 4-component link $(2,2,2)\,(2\,1,2\,1)$ with $n=12$
crossings (Fig. 1b) have the same Alexander polynomial $1-9 x+34
x^2-64 x^3+64 x^4-34 x^5+9 x^6-x^7$, and 3-component link
$6^*2\,2:.(2,-2)\,0$ with $n=12$ crossings and 5-component link
$2,2,2,2,2+$ with $n=11$ crossings have the same Alexander
polynomial $1-9 x+27 x^2-38 x^3+27 x^4-9 x^5+x^6$. However, up to
$n=12$ crossings the Alexander polynomial distinguishes links with
an odd number of components from links with even number of
components. Up to $n=12$ crossings all the remaining polynomials
completely distinguish links according to the number of components.

In the book [2], for families of alternating $KL$s we proposed the
following conjecture:

\begin{conjecture}
For every two alternating nonisotopic $KL$s $L_1$ and $L_2$
belonging to the same family $F$, $P(L_1)\neq P(L_2)$ for every
polynomial invariant $P$.
\end{conjecture}

From the obtained results we conclude that amount of all $KL$s with
$n\le 12$ crossings that cannot be detected by the mentioned
polynomial invariants is between $69\%$ (Alexander polynomial) and
$9\%$ (Kauffman two-variable polynomial). In this amount are included mutant
$KL$s that can not be distinguished by any polynomial invariant.

Comparing the results from Table 3 and Table 4 we conclude that for
all polynomials, except for the Kauffman polynomial the results are
worst if alternating and non-alternating $KL$s are not separated
before the computations, i.e., that for all polynomials, except for
the Kauffman polynomial there exist pairs (or groups) of $KL$s with
the same polynomial, which contain alternating and non-alternating
$KL$s. Up to $n=12$ crossings, every two $KL$s with the same
Kauffman polynomial have the same number of crossings. Hence, we
have the following open problem:

\medskip

\noindent {\bf Open problem 1:} {\it Find an alternating $KL$ with
the same Kauffman polynomial as some other non-alternating $KL$.}

\section{Factorizability of $KL$ polynomials and $KL$ \\ primeness detection}

The other test we made is the factorization, i.e., the ability of an
invariant to detect primeness of $KL$s. For all polynomial
invariants $P$, except the Khovanov polynomial
$P(L_1\#L_2)=P(L_1)P(L_2)$. However, the mentioned polynomials are
factorizable for some prime $KL$s as well. For example, the Jones
polynomial is factorizable for the link family 6, 10, $\ldots $,
$4k+2$, and for the rational knots $3\,1\,1\,3$ ($8_9$) (Fig. 2a),
$7\,2$ ($9_2$) (Fig. 2b), Homflypt polynomial is factorizable for
the 2-component link $2\,1\,1\,1\,1\,2$ ($8_8^2$) (Fig. 2c), and for knot $4\,2\,1\,2$ ($9_{12}$)
2-colored and 3-colored Jones polynomials are factorizable for the
3-component link $6,2,-2$ (Fig. 2d), {\it etc.} The only exceptions
we found are Tutte polynomial\footnote{Tutte polynomial is not $KL$
invariant, because it is not invariant under Reidemeister moves, but
it can be considered as the invariant of particular minimal
diagrams of alternating $KL$s.} and Kauffman two-variable
polynomial.

\begin{figure}
\centerline{\psfig{file=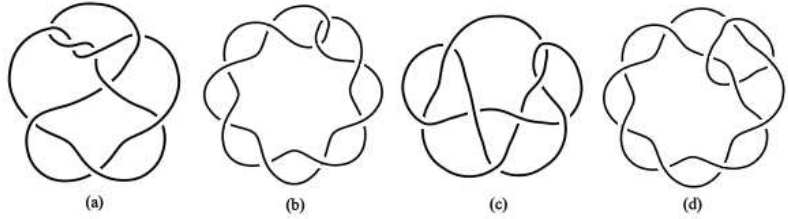,width=4.20in}} \vspace*{8pt}
\caption{(a) Knot $3\,1\,1\,3$; (b) knot $7\,2$; (c) 2-component
link $2\,1\,1\,1\,1\,2$; (d) 3-component link $6,2,-2$.} \label{f2}
\end{figure}

\begin{conjecture}
The Tutte polynomial and Kauffman two-variable polynomial detect
primeness, i.e, they are not factorizable for prime $KL$s.
\end{conjecture}

We expect that the conjecture about Tutte polynomial can be proved
on the basis of the irreducibility of the Tutte polynomial of
connected matroids (Brylawski theorem) [5]. Trying to find the
counterexample to the conjecture about Kauffman polynomial we
checked without success all rational $KL$s up to $n=19$ crossings,
all Montesinos $KL$s up to $n=18$ crossings, all knots from {\it
Knotscape} tables up to $n=16$ crossings, and all links up to $n=12$
crossings.


\begin{thebibliography}{5}

\bibitem{1}
Hoste, J. and Thistlethwaite, M. {\it Knotscape},
\newline \texttt{http://www.math.utk.edu/~morwen/knotscape.html}

\bibitem{2} S.V. Jablan and R. Sazdanovi\' c, {\it LinKnot -- Knot Theory by
Computer},  (World Scientific, New Jersey, London, Singapore, 2007;
\texttt{http://math.ict.edu.rs/},
\texttt{http://www.mi.sanu.ac.rs/vismath/linknot/index.html}).

\bibitem{3}
Shumakovitch, A. (2008) {\it KhoHo},
\texttt{http://www.geometrie.ch/KhoHo/}

\bibitem{4}
The Mathematica Package {\it KnotTheory},
\newline
\texttt{http://katlas.math.toronto.edu/wiki/The$_{-}$Mathematica$_{-}$Package$_{-}$}
\texttt{KnotTheory}

\bibitem{5} Merino, C., de Mier, A. and Noy, M. {\it Irreducibility of the Tutte polynomial of a connected matroid},
Journal of Combinatorial Theory, Series B, {\bf 83}, 2 (2001)
298--304.


\end{thebibliography}
\end{document}